\documentclass[10pt]{article}
\usepackage{amsmath}
\usepackage{amssymb}
\usepackage{amsfonts}
\usepackage{eucal}
\numberwithin{equation}{section} \oddsidemargin=0cm
\def\bincoeff#1#2{{#1\choose #2}}
\newtheorem{theorem}{Theorem}[section]

\newtheorem{lemma}[theorem]{Lemma}

\newtheorem{proposition}[theorem]{Proposition}
\newtheorem{remark}{Remark}[section]

\pdfoutput=1
\begin{document}
\title{On the asymptotic behavior of the hyperbolic Brownian motion}
\author{Valentina Cammarota\footnote{Dipartimento di Matematica,
Universit\`a degli Studi di Roma `Tor Vergata', Via della Ricerca
Scientifica 1, 00133 Roma, Italy. \textit{E-mail address}:
\texttt{cammarot@mat.uniroma2.it}, \textit{Web page:}
\texttt{https://sites.google.com/site/valentinacammarota/}} \and
Alessandro De Gregorio\footnote{Dipartimento di Scienze
Statistiche, Sapienza Universit\`{a} di Roma, Piazzale Aldo Moro
5, 00185 Rome, Italy. e-mail:
\texttt{alessandro.degregorio@uniroma1.it}} \and Claudio
Macci\footnote{Dipartimento di Matematica, Universit\`a degli
Studi di Roma `Tor Vergata', Via della Ricerca Scientifica 1,
00133 Roma, Italy. \textit{E-mail address}:
\texttt{macci@mat.uniroma2.it}}}
\date{}
\maketitle
\begin{abstract}
\noindent The main results in this paper concern large and
moderate deviations for the radial component of a $n$-dimensional
hyperbolic Brownian motion (for $n\geq 2$) on the Poincar\'{e}
half-space. We also investigate the asymptotic behavior of the
hitting probability $P_\eta(T_{\eta_1}^{(n)}<\infty)$ of a ball of
radius $\eta_1$, as the distance $\eta$ of the starting point of
the hyperbolic Brownian motion goes to infinity.\\
\ \\
\noindent\emph{Keywords}: hitting probability, hyperbolic
distance, large deviations, moderate deviations, Poincar\'{e}
half-space.\\
\noindent\emph{2000 Mathematical Subject Classification}: 60F10,
58J65, 60J60.
\end{abstract}

\section{Introduction}\label{sec:introduction}
Random motions in hyperbolic spaces, i.e. Riemannian manifolds
with constant negative curvature, have been studied since the
fifties and much attention has been placed on the so-called
hyperbolic Brownian motion on the Poincar\'e half-space
$\mathbb{H}^n$; the interested reader can consult, for example,
Gertsenshtein and Vasiliev \cite{GertsenshteinVasiliev}, Getoor
\cite{Getoor}, Gruet \cite{Gruet}, Matsumoto and Yor
\cite{MatsumotoYor}, Lao and Orsingher \cite{LaoOrsingher},
Byczkowski and Malecki \cite{ByczkowskiMalecki}, Borodin
\cite{Borodin}. Branching hyperbolic Brownian motion has been
analyzed by Lalley and Sellke \cite{Lalley} who investigated the
connection between the birth rate and the underlying dynamics in
supercritical and subcritical cases. Also Kelbert and Suhov
\cite{KelbertSuhov} and Karpelevich et al.
\cite{KarpelevichPecherskySuhov} have studied the asymptotic
behavior of the hyperbolic branching Brownian motion. Random walks
on the geodesics of the hyperbolic plane has been considered, for
example, in J\o rgensen \cite{Jorgensen} and Cammarota and
Orsingher \cite{CammarotaOrsingher2008}. One-dimensional and
planar random motions in non-Euclidean spaces have also been
analyzed in De Gregorio and Orsingher \cite{DegregorioOrsingher}.

There is a close link between one-dimensional disordered systems
and Brownian diffusion on hyperbolic spaces. For instance,
Gertsenshtein and Vasiliev \cite{GertsenshteinVasiliev}, in their
pioneering work, have shown that the statistical properties of
reflection and transmission coefficients for waveguides with
random inhomogeneities are directly related to some random walk on
the Poincar\'e half-plane.

Comtet and Monthus \cite{ComtetMonthus} showed how the
one-dimensional classical diffusion of a particle in a quenched
random potential is directly related to Brownian motion on the
hyperbolic plane. In particular the authors discussed some
functionals governing transport properties of a diffusion in a
random Brownian potential, in terms of hyperbolic Brownian motion.

The geodesic curves in Poincar\'e  half-plane model are either
half-circles with center lying on the $x$-axis or vertical
half-lines. For an optical non-homogenous media where light rays
move with velocity $c(x,y)=y$ (independent from direction), on the
base of Fermat's principle, the possible paths for the light are
those curves $L$ which satisfy the equality
\begin{eqnarray}\label{eq:1.1bis}
\frac{\sin \alpha (y)}{y}=k,
\end{eqnarray}
where $\alpha(y)$ is the angle between the vertical and the
tangent to $L$ in the point with ordinate $y$. It is easy to see
that the circles with center on the $x$-axis and radius
$\frac{1}{k}$ satisfy \eqref{eq:1.1bis} (for $k=0$ we get the
vertical lines). The hyperbolic geometrical optics where the ray
trajectories are the geodesics in the Poincar\'e  half-plane is
analyzed for example in De Micheli et al. \cite{DeMichelietal}.
Scattered obstacles in the non-homogeneous medium cause random
deviations in the propagation of light and this leads to the
random model analyzed in Cammarota and Orsingher
\cite{CammarotaOrsingher2008}.

Hyperbolic Brownian motion has been revitalized by mathematical
finance since some exotic derivatives have a strict connection
with the stochastic representation of the hyperbolic Brownian
motion. For example, Matsumoto and Yor \cite{Matsumoto} present
some identities for the pricing formula of the Asian or average
call option in the framework of the Black-Scholes model.

In this paper we present asymptotic results in the fashion of
large deviations. The theory of large deviations gives an
asymptotic computation of small probabilities on exponential scale
and it is used in several fields, and, in particular, in physics;
for instance a large number of authors (see e.g. Ellis
\cite{Ellis}) saw large deviation theory as the proper
mathematical framework in which problems of statistical mechanics
can be formulated and solved. We remark that, in analogy with what
happens in Monthus and Texier \cite{MonthusTexier} for random
walks on Bethe lattices, our asymptotic results can have interest
in the study of some random walks in the context of polymer
physics which admit the hyperbolic Brownian motion as a continuous
approximation.

We conclude with the outline of the paper. In Section
\ref{sec:preliminaries} we present some preliminaries on large
deviations and on the hyperbolic Brownian motion. In Section
\ref{sec:LDMD} we prove large and moderate deviation results for
the radial component of a $n$-dimensional hyperbolic Brownian
motion (for $n\geq 2$). In Section \ref{sec:hitting-probabilities}
we investigate the asymptotic behavior of the hitting
probabilities of an hyperbolic ball centered at the origin as the
distance of the starting point of the process goes to infinity.
Finally, in Section \ref{sec:Hirao}, we discuss some connections
with the results in Hirao \cite{Hirao}.

\section{Preliminaries}\label{sec:preliminaries}
We give some preliminaries on large deviations and on the
hyperbolic Brownian motion.

\subsection{Preliminaries on large deviations}
We recall the basic definitions in Dembo and Zeitouni
\cite{DemboZeitouni}, pages 4--5. Let $\mathcal{X}$ be a Hausdorff
topological space with Borel $\sigma$-algebra
$\mathcal{B}_{\mathcal{X}}$. A lower semi-continuous function
$I:\mathcal{X}\to [0,\infty]$ is called rate function. A family of
$\mathcal{X}$-valued random variables $\{X(t):t>0\}$ satisfies the
\emph{large deviation principle} (LDP for short), as $t\to\infty$,
with rate function $I$ and speed $v_t$ if:
$\lim_{t\to\infty}v_t=\infty$,
$$\limsup_{t\to\infty}\frac{1}{v_t}\log P(X(t)\in F)\leq-\inf_{x\in
F}I(x)\quad\textrm{for all closed sets}\ F$$ and
$$\liminf_{t\to\infty}\frac{1}{v_t}\log P(X(t)\in G)\geq-\inf_{x\in
G}I(x)\quad\textrm{for all open sets}\ G.$$ A rate function $I$ is
said to be good if all the level sets
$\{\{x\in\mathcal{X}:I(x)\leq\gamma\}:\gamma\geq 0\}$ are compact.
We recall that the lower bound for open sets is equivalent to the
following condition (see eq. (1.2.8) in Dembo and Zeitouni
\cite{DemboZeitouni}):
\begin{equation}\label{eq:LB-local-condition}
\liminf_{t\to\infty}\frac{1}{v_t}\log P(X(t)\in G)\geq-I(x)\quad
\left.\begin{array}{ll}
\mbox{for all $x\in\mathcal{X}$ such that $I(x)<\infty$ and}\\
\mbox{for all open sets $G$ such that $x\in G$}.
\end{array}\right.
\end{equation}
In this paper we prove LDPs with $\mathcal{X}=\mathbb{R}$. We
start with Proposition \ref{prop:LDP} where we have $v_t=t$. We
also study the \emph{moderate deviations}, i.e. we prove
Proposition \ref{prop:MDP} which provides a class of LDPs where
$v_t=t^{1-2\beta}$, varying $\beta\in(0,1/2)$, and the rate
function does not depend on $\beta$. In some sense the moderate
deviations fill the gap between an asymptotic normality result
(i.e. Theorem 2.1 in Matsumoto \cite{Matsumoto}) for
$\beta=\frac{1}{2}$, and a convergence in probability to a
constant (i.e. Corollary 5.7.3 in Davies \cite{Davies}) - together
with a centering of the random variables - for $\beta=0$. We
remark that we have a quadratic rate function which vanishes at
zero. To better explain this concept, we recall the basic result
on moderate deviations for the empirical means
$\left\{\frac{X_1+\cdots+X_n}{n}:n\geq 1\right\}$ of i.i.d.
centered and $\mathbb{R}^d$-valued random variables $\{X_n:n\geq
1\}$ (see e.g. Theorem 3.7.1 in Dembo and Zeitouni
\cite{DemboZeitouni}) which fill the gap between the central limit
theorem and the law of the large numbers. In such a case we have
the LDP for $\left\{\sqrt{na_n}\frac{X_1+\cdots+X_n}{n}:n\geq
1\right\}$ for $\{a_n:n\geq 1\}$ such that $a_n\to 0$ and
$na_n\to\infty$ (as $n\to\infty$) with speed $v_n=\frac{1}{a_n}$;
then one can check that $a_n$ plays the role of $t^{2\beta-1}$ in
Proposition \ref{prop:MDP} (actually $t^{2\beta-1}\to 0$ and
$t\cdot t^{2\beta-1}\to\infty$ as $t\to\infty$).

\subsection{Preliminaries on the hyperbolic Brownian motion}
For $n\geq 2$, let
$\mathbb{H}^n=\{z=(x,y):x\in\mathbb{R}^{n-1},y>0\}$ be the upper
half-space with origin $O_n=(0,\ldots,0,1)$ endowed with the
hyperbolic Riemannian metric
$$ds^2=\frac{d x_1^2+\cdots+d x_{n-1}^2+d y^2}{y^2}.$$
The hyperbolic distance $\mathrm{\eta}(z,z')$ between $z=(x,y)$
and $z'=(x',y')$ in $\mathbb{H}^n$ is given by the formula
$$\cosh \mathrm{\eta}(z,z')=\frac{|x-x'|^2+y^2+y'^2}{2yy'}$$
where $|x-x'|$ is the Euclidean distance between $x,x' \in
\mathbb{R}^{n-1}$, and the volume element is given by
$$dv={y^{-n}} \,{dx dy}=\sinh^{n-1} \eta \; d \eta \; d\Omega_n$$
where $d\Omega_n$ is the surface element of the $n$-dimensional
unit sphere. The Laplace-Beltrami operator in $\mathbb{H}^n$ is
\begin{equation*}
\Delta_n=y^2 \left( \sum_{i=1}^{n-1} \frac{\partial^2}{\partial x_i^2}+\frac{\partial^2}{\partial y^2}\right)-\frac{(n-2)}2 y \frac{\partial}{\partial y}
\end{equation*}
(see, for example, Chavel \cite{Chavel} page 265).

The hyperbolic Brownian motion is a diffusion governed by the
generator ${\Delta_n}/2$. We denote by $k_n(z,z',t)$ the heat
kernel, with respect to the volume element $dv$. Since the
Laplace-Beltrami operator is invariant  under diffeomorphism,
$k_n(z,z',t)$ is a function of $\mathrm{\eta}(z,z')$ and we write
$k_n(\eta,t)$ for $k_n(z,z',t)$. Therefore the transition density
$p_n(\eta,t)$ of the radial component of hyperbolic Brownian
motion is given by
$$p_n(\eta,t)=\frac{2 \pi^{n/2}}{\Gamma(n/2)} k_n(\eta,t) \sinh^{n-1} \eta \; d \eta\ (\eta>0,\ t>0).$$
The classical formulae for the heat kernel, which are of different
forms for odd and even dimensions $n$, are well known together
with the recurrence Millson's formula with respect to the
dimension $n$ (see e.g Davies and Mandouvalos
\cite{DaviesMandouvalos}). The analogue formulae for the
hyperbolic heat kernel can be found in several references: see,
e.g., Buser \cite{Buser} Theorem 7.4.1, Chavel \cite{Chavel}
Section X.2, Terras \cite{Terras} Section 3.2 and Helgason
\cite{Helgason} page 29. In detail we have
\begin{align*}
k_2(\eta ,t)&=\frac{e^{-\frac{t}{4}}}{ 2^{5/2} ( \pi t)^{3/2} } \int_{\eta }^{\infty }\frac{\varphi e^{-\frac{\varphi ^{2}}{4t}}}{\sqrt{
\cosh \varphi -\cosh \eta }}d\varphi, \\
k_3(\eta ,t)&=\frac{e^{-t} }{2^3 (\pi t)^{3/2}} \frac{\eta e^{-
\frac{\eta^2}{4 t}}}{\sinh \eta}
\end{align*}
and, in general, for all $n\geq 2$, closed form expressions are
not available for $k_n=k_n(\eta,t)$. Therefore we will use some
known sharp bounds, see e.g. Theorem 5.7.2 in Davies \cite{Davies}
or Theorem 3.1 in Davies and Mandouvalos \cite{DaviesMandouvalos}
(note that $n$ in that reference is replaced by $n-1$ in this
paper). Let $h_n$ be defined by
$$h_n(\eta,t):=t^{-n/2}\exp\left(-\frac{(n-1)^2t}{4}-\frac{(n-1)\eta}{2}-\frac{\eta^2}{4t}\right)(1+\eta+t)^{(n-3)/2}(1+\eta);$$
we have that $k_n(\eta,t)\sim h_n(\eta)$, i.e. there exists
$c_n\in(1,\infty)$ such that we have
$$c_n^{-1}h_n(\eta,t)\leq k_n(\eta,t)\leq c_nh_n(\eta,t)$$
for all $t>0$ and $\eta>0$. Then, since
$(\sinh\eta)^{n-1}\sim\left(\frac{\eta}{1+\eta}\right)^{n-1}e^{(n-1)\eta}$,
if we set
$$g_n(\eta,t):=t^{-n/2}\left(\frac{\eta}{1+\eta}\right)^{n-1}
\exp\left(-\frac{(\eta-(n-1)t)^2}{4t}\right)(1+\eta+t)^{(n-3)/2}(1+\eta),$$
we can say that
$$(\sinh\eta)^{n-1}k_n(\eta,t)\sim g_n(\eta,t)$$
uniformly in $\eta$ and $t$, i.e. there exists a constant
$d_n\in(1,\infty)$ such that
\begin{equation}\label{eq:DM-consequence}
d_n^{-1}g_n(\eta,t)\leq(\sinh\eta)^{n-1}k_n(\eta,t)\leq
d_ng_n(\eta,t)
\end{equation}
for all $t>0$ and $\eta>0$.

\section{Large and moderate deviations}\label{sec:LDMD}
In this section we prove asymptotic results for the radial
component of a $n$-dimensional hyperbolic Brownian motion (for
$n\geq 2$), which will be denoted by $\{D_n(t):t\geq 0\}$. More
precisely we prove two LDPs: the first one (Proposition
\ref{prop:LDP}) concerns the convergence in probability to a
constant, the second one (Proposition \ref{prop:MDP}) concerns the
moderate deviation regime. Both proofs are divided in two parts.
\begin{enumerate}
\item The proof of the lower bound for open sets, and we refer
to condition \eqref{eq:LB-local-condition} with appropriate
choices of $\{X(t):t>0\}$, $v_t$ and $I$.
\item The proof of the upper bound for closed sets, and we
often refer to an upper bound for the moment generating function
$\mathbb{E}[e^{\lambda D_n(t)}]$ for all $\lambda\in\mathbb{R}$
(actually we have to consider $\lambda\neq 0$) which is given in
the next Lemma \ref{lemma:mgf}.
\end{enumerate}

\begin{lemma}\label{lemma:mgf}
Let $d_n$ be as in \eqref{eq:DM-consequence} and let
$\lambda\in\mathbb{R}$ be arbitrarily fixed. Moreover we set
$\kappa(\lambda):=\lambda(\lambda+n-1)$ and
$m_n:=\lceil\frac{n-1}{2}\rceil$. Then, for all $t>0$, we have
$$\mathbb{E}[e^{\lambda D_n(t)}]\leq
2\sqrt{2\pi(1/(2t))}d_nt^{-n/2+1}e^{\kappa(\lambda)t}\sum_{j=0}^{m_n}\bincoeff{m_n}{j}(1+nt)^{m_n-j}(2t)^j\mathbb{E}[|W_t|^j],$$
where $W_t$ is a Normal distributed random variable with mean
$\lambda$ and variance $\frac{1}{2t}$.
\end{lemma}
\noindent\emph{Proof.} Firstly, by \eqref{eq:DM-consequence}, we
have
\begin{align*}
\mathbb{E}[e^{\lambda D_n(t)}]\leq&d_n\int_0^\infty
e^{\lambda\eta}t^{-n/2}\left(\frac{\eta}{1+\eta}\right)^{n-1}\exp\left(-\frac{(\eta-(n-1)t)^2}{4t}\right)(1+\eta+t)^{(n-3)/2}(1+\eta)d\eta\\
\leq&d_nt^{-n/2}\int_0^\infty
e^{\lambda\eta}\exp\left(-\frac{(\eta-(n-1)t)^2}{4t}\right)(1+\eta+t)^{(n-1)/2}d\eta;
\end{align*}
thus, by the change of variable $\alpha=\frac{\eta-(n-1)t}{2t}$
and some computations, we obtain
\begin{align*}
\mathbb{E}[e^{\lambda
D_n(t)}]\leq&d_nt^{-n/2}e^{\kappa(\lambda)t}\int_0^\infty\exp\left(\lambda\eta-\lambda(\lambda+n-1)t-\frac{(\eta-(n-1)t)^2}{4t}\right)(1+\eta+t)^{(n-1)/2}d\eta\\
=&2d_nt^{-n/2+1}e^{\kappa(\lambda)t}A(\lambda,t,n),
\end{align*}
where
$$A(\lambda,t,n):=\int_{(1-n)/2}^\infty\exp\left(-(\lambda-\alpha)^2t\right)(1+2\alpha t+nt)^{(n-1)/2}d\alpha.$$
Finally, since
\begin{align*}
A(\lambda,t,n)\leq&\int_{(1-n)/2}^\infty\exp\left(-(\lambda-\alpha)^2t\right)(1+2\alpha
t+nt)^{m_n}d\alpha\\
=&\sum_{j=0}^{m_n}\bincoeff{m_n}{j}(1+nt)^{m_n-j}(2t)^j\int_{(1-n)/2}^\infty\exp\left(-(\lambda-\alpha)^2t\right)\alpha^jd\alpha\\
\leq&\sum_{j=0}^{m_n}\bincoeff{m_n}{j}(1+nt)^{m_n-j}(2t)^j\int_{-\infty}^\infty\exp\left(-(\lambda-\alpha)^2t\right)|\alpha|^jd\alpha
\end{align*}
and
$$\int_{-\infty}^\infty\exp\left(-(\lambda-\alpha)^2t\right)|\alpha|^jd\alpha
=\int_{-\infty}^\infty\exp\left(-\frac{(\lambda-\alpha)^2}{2(1/(2t))}\right)|\alpha|^jd\alpha=\sqrt{2\pi(1/(2t))}\mathbb{E}[|W_t|^j],$$
we conclude by putting together the inequalities for
$\mathbb{E}[e^{\lambda D_n(t)}]$ and $A(\lambda,t,n)$. $\Box$\\

We start with the LDP associated to the convergence in probability
of $\frac{D_n(t)}{t}$ to $n-1$ (as $t\to\infty$).

\begin{proposition}\label{prop:LDP}
The family $\left\{\frac{D_n(t)}{t}:t>0\right\}$ satisfies the LDP
with good rate function $I_1$ defined by
$$I_1(x)=\left\{\begin{array}{ll}
\frac{(x-(n-1))^2}{4}&\ \mbox{if}\ x\geq 0\\
\infty&\ \mbox{if}\ x<0,
\end{array}\right.$$
and speed function $v_t=t$.
\end{proposition}
\noindent\emph{Proof.} The proof is divided in two parts.\\
\emph{1) Lower bound for open sets.} We have to check that
$$\liminf_{t\to\infty}\frac{1}{t}\log P\left(\frac{D_n(t)}{t}\in G\right)\geq-\frac{(x-(n-1))^2}{4}$$
for all $x\geq 0$ and for all open sets $G$ such that $x\in G$. We
have two cases.\\
\noindent\emph{Case $x>0$.} Let $\varepsilon>0$ be such that
$(x-\varepsilon,x+\varepsilon)\subset G$; moreover we can choose
$\varepsilon\in(0,x)$. Then, by \eqref{eq:DM-consequence} and by
considering the change of variable $\alpha=\frac{\eta}{t}$, we
have
\begin{align*}
P\left(\frac{D_n(t)}{t}\in G\right)\geq&P\left(\frac{D_n(t)}{t}\in (x-\varepsilon,x+\varepsilon)\right)\\
\geq&d_n^{-1}\int_{(x-\varepsilon)t}^{(x+\varepsilon)t}t^{-n/2}\left(\frac{\eta}{1+\eta}\right)^{n-1}
\exp\left(-\frac{(\eta-(n-1)t)^2}{4t}\right)(1+\eta+t)^{(n-3)/2}\\
&\cdot(1+\eta)d\eta\\
\geq&d_n^{-1}t^{-n/2+1}\int_{x-\varepsilon}^{x+\varepsilon}\left(\frac{\alpha
t}{1+\alpha t}\right)^{n-1}
\exp\left(-\frac{(\alpha-(n-1))^2t}{4}\right)(1+\alpha
t+t)^{(n-3)/2}d\alpha\\
\geq&d_n^{-1}t^{-n/2+1}2\varepsilon\left(\frac{(x-\varepsilon)t}{1+(x+\varepsilon)t}\right)^{n-1}\\
&\cdot\exp\left(\min\left\{-\frac{(x+\varepsilon-(n-1))^2t}{4},-\frac{(x-\varepsilon-(n-1))^2t}{4}\right\}\right)\\
&\cdot(1+(x-\varepsilon)t+t)^{(n-3)/2};
\end{align*}
thus
$$\liminf_{t\to\infty}\frac{1}{t}\log P\left(\frac{D_n(t)}{t}\in G\right)
\geq-\max\left\{\frac{(x-\varepsilon-(n-1))^2}{4},\frac{(x+\varepsilon-(n-1))^2}{4}\right\},$$
and we conclude by letting $\varepsilon$ go to zero.\\
\noindent\emph{Case $x=0$.} Let $\{x_k:k\geq 1\}$ and
$\{\delta_k:k\geq 1\}$ be two sequences of positive numbers such
that $\lim_{k\to\infty}x_k=0$ and, for all $k\geq 1$,
$(x_k-\delta_k,x_k+\delta_k)\subset G$. Then, for each fixed
$k\geq 1$, we have
$$\liminf_{t\to\infty}\frac{1}{t}\log P\left(\frac{D_n(t)}{t}\in G\right)
\geq\liminf_{t\to\infty}\frac{1}{t}\log P\left(\frac{D_n(t)}{t}\in
(x_k-\delta_k,x_k+\delta_k)\right)\geq-\frac{(x_k-(n-1))^2}{4}$$
by the first part of the proof concerning $x>0$ (we have $x_k$ in
place of $x$ and $(x_k-\delta_k,x_k+\delta_k)$ in place of $G$)
and we complete the proof (for the case $x=0$) letting $k$ go to
infinity.\\
\emph{2) Upper bound for closed sets.} The upper bound
$$\limsup_{t\to\infty}\frac{1}{t}\log P\left(\frac{D_n(t)}{t}\in F\right)\leq-\inf_{x\in
F}I_1(x)\quad\textrm{for all closed sets}\ F$$ trivially holds if
$n-1\in F$ or $F\cap [0,\infty)$ is empty. Thus, from now on, we
assume that $n-1\notin F$ and $F\cap [0,\infty)$ is nonempty. We
also assume that both $F\cap [0,n-1)$ and $F\cap (n-1,\infty)$ are
nonempty; actually at least one of the two sets is nonempty and,
if one of them would be empty, the proof presented below could be
readily adapted. We define
$$\hat{x}:=\sup(F\cap [0,n-1))\ \textrm{and}\ \tilde{x}:=\inf(F\cap (n-1,\infty))$$
and we have $\hat{x},\tilde{x}\in F$, $0\leq
\hat{x}<n-1<\tilde{x}$, and $F\subset
(-\infty,\hat{x}]\cup[\tilde{x},\infty)$. Then
$$P\left(\frac{D_n(t)}{t}\in F\right)\leq P\left(\frac{D_n(t)}{t}\leq \hat{x}\right)+P\left(\frac{D_n(t)}{t}\geq \tilde{x}\right)$$
and, by Lemma 1.2.15 in Dembo and Zeitouni \cite{DemboZeitouni},
we get
\begin{align*}
\limsup_{t\to\infty}\frac{1}{t}\log
&P\left(\frac{D_n(t)}{t}\in F\right)\\
\leq&\max\left\{\limsup_{t\to\infty}\frac{1}{t}\log
P\left(\frac{D_n(t)}{t}\leq
\hat{x}\right),\limsup_{t\to\infty}\frac{1}{t}\log
P\left(\frac{D_n(t)}{t}\geq \tilde{x}\right)\right\}.
\end{align*}
Thus, if we prove
\begin{equation}\label{eq:UB1LDP}
\limsup_{t\to\infty}\frac{1}{t}\log P\left(\frac{D_n(t)}{t}\leq
\hat{x}\right)\leq-\inf_{x\leq
\hat{x}}I_1(x)=-\frac{(\hat{x}-(n-1))^2}{4}
\end{equation}
and
\begin{equation}\label{eq:UB2LDP}
\limsup_{t\to\infty}\frac{1}{t}\log P\left(\frac{D_n(t)}{t}\geq
\tilde{x}\right)\leq-\inf_{x\geq
\tilde{x}}I_1(x)=-\frac{(\tilde{x}-(n-1))^2}{4},
\end{equation}
we conclude the proof because we get
\begin{align*}
\limsup_{t\to\infty}\frac{1}{t}\log P\left(\frac{D_n(t)}{t}\in
F\right)\leq&\max\left\{-\inf_{x\leq \hat{x}}I_1(x),-\inf_{x\geq
\tilde{x}}I_1(x)\right\}\\
\leq&-\min\left\{\inf_{x\leq \hat{x}}I_1(x),\inf_{x\geq
\tilde{x}}I_1(x)\right\}=-\inf_{x\in F}I_1(x).
\end{align*}

\noindent\emph{Proof of \eqref{eq:UB1LDP}.} The case $\hat{x}=0$
is trivial because we have $-\infty\leq-\frac{(n-1)^2}{4}$ noting
that $P\left(\frac{D_n(t)}{t}\leq 0\right)=0$. For $\hat{x}>0$, by
\eqref{eq:DM-consequence}, we have
\begin{align*}
P\left(\frac{D_n(t)}{t}\leq
\hat{x}\right)\leq&d_n\int_0^{t\hat{x}}t^{-n/2}\left(\frac{\eta}{1+\eta}\right)^{n-1}
\exp\left(-\frac{(\eta-(n-1)t)^2}{4t}\right)(1+\eta+t)^{(n-3)/2}(1+\eta)d\eta\\
\leq&d_nt^{-n/2}t\hat{x}\exp\left(-\frac{(\hat{x}-(n-1))^2t}{4}\right)(1+t\hat{x}+t)^{(n-3)/2}(1+t\hat{x}),
\end{align*}
and we obtain
$$\limsup_{t\to\infty}\frac{1}{t}\log P\left(\frac{D_n(t)}{t}\leq\hat{x}\right)\leq-\frac{(\hat{x}-(n-1))^2}{4}.$$
\noindent\emph{Proof of \eqref{eq:UB2LDP}.} Firstly, by Markov's
inequality, for all $\lambda>0$ we have
$$P\left(\frac{D_n(t)}{t}\geq\tilde{x}\right)\leq\mathbb{E}[e^{\lambda D_n(t)}]e^{-\lambda t\tilde{x}}.$$
Moreover, by Lemma \ref{lemma:mgf}, we get
$$P\left(\frac{D_n(t)}{t}\geq
\tilde{x}\right)\leq e^{-\lambda
t\tilde{x}}2\sqrt{2\pi(1/(2t))}d_nt^{-n/2+1}e^{\kappa(\lambda)t}\sum_{j=0}^{m_n}\bincoeff{m_n}{j}(1+nt)^{m_n-j}(2t)^j\mathbb{E}[|W_t|^j].$$
Now let $t_0>0$ be arbitrarily fixed; then, for all $t>t_0$, we
have
\begin{align*}
\mathbb{E}[|W_t|^j]=&\frac{1}{\sqrt{2\pi(1/(2t))}}\int_{-\infty}^\infty\exp\left(-\frac{(\lambda-\alpha)^2}{2(1/(2t))}\right)|\alpha|^jd\alpha\\
\leq&\sqrt{\frac{t}{t_0}}\frac{1}{\sqrt{2\pi(1/(2t_0))}}\int_{-\infty}^\infty\exp\left(-\frac{(\lambda-\alpha)^2}{2(1/(2t_0))}\right)|\alpha|^jd\alpha
=\sqrt{\frac{t}{t_0}}\mathbb{E}[|W_{t_0}|^j].
\end{align*}
In conclusion we have
$$\limsup_{t\to\infty}\frac{1}{t}\log P\left(\frac{D_n(t)}{t}\geq \tilde{x}\right)\leq-\lambda \tilde{x}+\kappa(\lambda),$$
and therefore
$$\limsup_{t\to\infty}\frac{1}{t}\log P\left(\frac{D_n(t)}{t}\geq
\tilde{x}\right)\leq\inf_{\lambda>0}\{-\lambda
\tilde{x}+\kappa(\lambda)\}=-\inf_{x\geq \tilde{x}}I_1(x)$$
because $\inf_{\lambda>0}\{-\lambda
\tilde{x}+\kappa(\lambda)\}=-\frac{(\tilde{x}-(n-1))^2}{4}=-I_1(\tilde{x})$
(actually the infimum is attained at
$\lambda=\frac{\tilde{x}-(n-1)}{2}$). $\Box$

\begin{remark}\label{rem:analogue-Euclidean-result}
We can state the analogue of Proposition \ref{prop:LDP} for a
centered Euclidean $n$-dimensional Brownian motion $\{B_n(t):t\geq
0\}$. One can easily check that
$\left\{\frac{\|B_n(t)\|}{t}:t>0\right\}$ satisfies the LDP with
good rate function $J_1$ defined by
$$J_1(x):=\left\{\begin{array}{ll}
\frac{x^2}{2}&\ \mbox{if}\ x\geq 0\\
\infty&\ \mbox{if}\ x<0
\end{array}\right.$$
with a standard application of the G\"{a}rtner Ellis Theorem (see
e.g. Theorem 2.3.6 in Dembo and Zeitouni \cite{DemboZeitouni}) and
of the contraction principle (see e.g. Theorem 4.2.1 in Dembo and
Zeitouni \cite{DemboZeitouni}). Thus, in some sense (see also the
note just after the statement of Corollary 5.7.3 in Davies
\cite{Davies}), the hyperbolic Brownian motion has an implicit
drift directed away from the origin because $I_1(x)$ is a
quadratic rate function (on $[0,\infty)$) which vanishes at
$x=n-1$, while $J_1(x)$ vanishes at $x=0$.
\end{remark}

We conclude with the LDPs concerning the moderate deviation
regime.

\begin{proposition}\label{prop:MDP}
For all $\beta\in(0,1/2)$, the family
$\left\{t^{(2\beta-1)/2}\left(\frac{D_n(t)-(n-1)t}{\sqrt{t}}\right):t>0\right\}$
satisfies the LDP with rate function $I_2$ defined by
$I_2(x)=\frac{x^2}{4}$ (for $x\in\mathbb{R}$) and speed function
$v_t=t^{1-2\beta}$.
\end{proposition}
\noindent\emph{Proof.} Firstly, in order to have simpler formulae,
we remark that
$$t^{(2\beta-1)/2}\frac{D_n(t)-(n-1)t}{\sqrt{t}}=t^{\beta-1}(D_n(t)-(n-1)t).$$
The proof is divided in two parts.\\
\emph{1) Lower bound for open sets.} We have to check that
$$\liminf_{t\to\infty}\frac{1}{t^{1-2\beta}}\log P\left(t^{\beta-1}(D_n(t)-(n-1)t)\in G\right)\geq-\frac{x^2}{4}$$
for all $x\in\mathbb{R}$ and for all open sets $G$ such that $x\in
G$. We can find $\varepsilon>0$ such that
$(x-\varepsilon,x+\varepsilon)\subset G$. Then we have
\begin{align*}
P\left(t^{\beta-1}(D_n(t)-(n-1)t)\in G\right)
\geq&P\left(t^{\beta-1}(D_n(t)-(n-1)t)\in (x-\varepsilon,x+\varepsilon)\right)\\
\geq&P\left(D_n(t)\in\left(\frac{x-\varepsilon}{t^{\beta-1}}+(n-1)t,\frac{x+\varepsilon}{t^{\beta-1}}+(n-1)t\right)\right)
\end{align*}
and, from now on, we take $t$ large enough to have
$\frac{x-\varepsilon}{t^{\beta-1}}+(n-1)t>0$. Now, by
\eqref{eq:DM-consequence} and by considering the change of
variable $\alpha=t^{\beta-1}(\eta-t(n-1))$, we have
\begin{align*}
P\left(t^{\beta-1}(D_n(t)-(n-1)t)\in G\right)
\geq&d_n^{-1}\int_{\frac{x-\varepsilon}{t^{\beta-1}}+(n-1)t}^{\frac{x+\varepsilon}{t^{\beta-1}}+(n-1)t}t^{-n/2}\left(\frac{\eta}{1+\eta}\right)^{n-1}\\
&\cdot\exp\left(-\frac{(\eta-(n-1)t)^2}{4t}\right)(1+\eta+t)^{(n-3)/2}(1+\eta)d\eta\\
\geq&d_n^{-1}t^{-n/2+1-\beta}\int_{x-\varepsilon}^{x+\varepsilon}\left(\frac{\frac{\alpha}{t^{\beta-1}}+tn-t}{1+\frac{\alpha}{t^{\beta-1}}+tn-t}\right)^{n-1}\\
&\cdot\exp\left(-\frac{\alpha^2}{4t^{2\beta-1}}\right)\left(1+\frac{\alpha}{t^{\beta-1}}+tn\right)^{(n-3)/2}d\alpha.
\end{align*}
Moreover there exists
$\alpha_{\varepsilon,t}\in(x-\varepsilon,x+\varepsilon)$ such that
\begin{align*}
P\left(t^{\beta-1}(D_n(t)-(n-1)t)\in G\right)
\geq&d_n^{-1}t^{-n/2+1-\beta}2\varepsilon\left(\frac{\frac{\alpha_{\varepsilon,t}}{t^{\beta-1}}+tn-t}
{1+\frac{\alpha_{\varepsilon,t}}{t^{\beta-1}}+tn-t}\right)^{n-1}\\
&\cdot\exp\left(-\frac{\alpha_{\varepsilon,t}^2}{4t^{2\beta-1}}\right)\left(1+\frac{\alpha_{\varepsilon,t}}{t^{\beta-1}}+tn\right)^{(n-3)/2};
\end{align*}
this yields
$$\liminf_{t\to\infty}\frac{1}{t^{1-2\beta}}\log P\left(t^{\beta-1}(D_n(t)-(n-1)t)\in
G\right)\geq-\max\left\{\frac{(x-\varepsilon)^2}{4},\frac{(x+\varepsilon)^2}{4}\right\},$$
and we conclude by letting $\varepsilon$ go to zero.\\
\emph{2) Upper bound for closed sets.} The upper bound
$$\limsup_{t\to\infty}\frac{1}{t^{1-2\beta}}\log P\left(t^{\beta-1}(D_n(t)-(n-1)t)\in F\right)\leq-\inf_{x\in
F}I_2(x)\quad\textrm{for all closed sets}\ F$$ trivially holds if
$0\in F$ or $F$ is empty. Thus, from now on, we assume that
$0\notin F$ and $F$ is nonempty. We also assume that both $F\cap
(-\infty,0)$ and $F\cap (0,\infty)$ are nonempty; actually at
least one of the two sets is nonempty and, if one of them would be
empty, the proof presented below could be readily adapted. We
define
$$\hat{x}:=\sup(F\cap (-\infty,0))\ \textrm{and}\ \tilde{x}:=\inf(F\cap (0,\infty))$$ and we
have $\hat{x},\tilde{x}\in F$, $\hat{x}<0<\tilde{x}$, and
$F\subset (-\infty,\hat{x}]\cup[\tilde{x},\infty)$. Then
\begin{align*}
P\left(t^{\beta-1}(D_n(t)-(n-1)t)\in F\right)\leq&P\left(t^{\beta-1}(D_n(t)-(n-1)t)\leq\hat{x}\right)\\
&+P\left(t^{\beta-1}(D_n(t)-(n-1)t)\geq \tilde{x}\right)
\end{align*}
and, by Lemma 1.2.15 in Dembo and Zeitouni \cite{DemboZeitouni},
we get
\begin{align*}
\limsup_{t\to\infty}\frac{1}{t^{1-2\beta}}\log&
P\left(t^{\beta-1}(D_n(t)-(n-1)t)\in
F\right)\\
\leq&\max\left\{\limsup_{t\to\infty}\frac{1}{t^{1-2\beta}}\log
P\left(t^{\beta-1}(D_n(t)-(n-1)t)\leq\hat{x}\right),\right.\\
&\left.\limsup_{t\to\infty}\frac{1}{t^{1-2\beta}}\log
P\left(t^{\beta-1}(D_n(t)-(n-1)t)\geq \tilde{x}\right)\right\}.
\end{align*}
Thus, if we prove
\begin{equation}\label{eq:UB1MDP}
\limsup_{t\to\infty}\frac{1}{t^{1-2\beta}}\log
P\left(t^{\beta-1}(D_n(t)-(n-1)t)\leq
\hat{x}\right)\leq-\inf_{x\leq \hat{x}}I_2(x)=-\frac{\hat{x}^2}{4}
\end{equation}
and
\begin{equation}\label{eq:UB2MDP}
\limsup_{t\to\infty}\frac{1}{t^{1-2\beta}}\log
P\left(t^{\beta-1}(D_n(t)-(n-1)t)\geq
\tilde{x}\right)\leq-\inf_{x\geq
\tilde{x}}I_2(x)=-\frac{\tilde{x}^2}{4},
\end{equation}
we conclude the proof because we get
\begin{align*}
\limsup_{t\to\infty}\frac{1}{t^{1-2\beta}}\log
P\left(t^{\beta-1}(D_n(t)-(n-1)t)\in
F\right)\leq&\max\left\{-\inf_{x\leq \hat{x}}I_2(x),-\inf_{x\geq
\tilde{x}}I_2(x)\right\}\\
\leq&-\min\left\{\inf_{x\leq \hat{x}}I_2(x),\inf_{x\geq
\tilde{x}}I_2(x)\right\}=-\inf_{x\in F}I_2(x).
\end{align*}
\noindent\emph{Proof of \eqref{eq:UB1MDP}.} Firstly, by Markov's
inequality, for all $\lambda<0$ we have
$$P\left(t^{\beta-1}(D_n(t)-(n-1)t)\leq\hat{x}\right)=P(D_n(t)\leq\hat{x}t^{1-\beta}+t(n-1))\leq\mathbb{E}[e^{\lambda
D_n(t)}]e^{-\lambda(\hat{x}t^{1-\beta}+t(n-1))}.$$ Moreover, by
Lemma \ref{lemma:mgf}, we get
\begin{align*}
P\left(t^{\beta-1}(D_n(t)-(n-1)t)\leq \hat{x}\right)\leq&
e^{-\lambda(\hat{x}t^{1-\beta}+t(n-1))}2\sqrt{2\pi(1/(2t))}d_nt^{-n/2+1}e^{\kappa(\lambda)t}\\
&\cdot\sum_{j=0}^{m_n}\bincoeff{m_n}{j}(1+nt)^{m_n-j}(2t)^j\mathbb{E}[|W_t|^j].
\end{align*}
We remark that, if we argue as in the proof of \eqref{eq:UB2LDP}
(in the proof of Proposition \ref{prop:LDP}), we obtain the
estimate
$\mathbb{E}[|W_t|^j]\leq\sqrt{\frac{t}{t_0}}\mathbb{E}[|W_{t_0}|^j]$
for some arbitrarily fixed $t_0>0$ and for $t>t_0$. Finally we
take $\lambda=\frac{\hat{x}}{2}t^{-\beta}$, and we have
\begin{align*}
\limsup_{t\to\infty}\frac{1}{t^{1-2\beta}}\log
&P\left(t^{\beta-1}(D_n(t)-(n-1)t)\leq
\hat{x}\right)\\
\leq&\limsup_{t\to\infty}\frac{1}{t^{1-2\beta}}\left\{-\frac{\hat{x}}{2}t^{-\beta}(\hat{x}t^{1-\beta}+t(n-1))+
\kappa\left(\frac{\hat{x}}{2}t^{-\beta}\right)t\right\}\\
=&\limsup_{t\to\infty}\left\{-\frac{\hat{x}^2}{2}-\frac{\hat{x}}{2}(n-1)t^\beta+
\frac{\hat{x}}{2}t^\beta\left(\frac{\hat{x}}{2}t^{-\beta}+n-1\right)\right\}=-\frac{\hat{x}^2}{2}+\frac{\hat{x}^2}{4}=-\frac{\hat{x}^2}{4}.
\end{align*}
\noindent\emph{Proof of \eqref{eq:UB2MDP}.} Firstly, by Markov's
inequality, for all $\lambda>0$ we have
$$P\left(t^{\beta-1}(D_n(t)-(n-1)t)\geq\tilde{x}\right)=P(D_n(t)\geq\tilde{x}t^{1-\beta}+t(n-1))\leq\mathbb{E}[e^{\lambda
D_n(t)}]e^{-\lambda(\tilde{x}t^{1-\beta}+t(n-1))}.$$ Moreover, by
Lemma \ref{lemma:mgf}, we get
\begin{align*}
P\left(t^{\beta-1}(D_n(t)-(n-1)t)\geq \tilde{x}\right)\leq&
e^{-\lambda(\tilde{x}t^{1-\beta}+t(n-1))}2\sqrt{2\pi(1/(2t))}d_nt^{-n/2+1}e^{\kappa(\lambda)t}\\
&\cdot\sum_{j=0}^{m_n}\bincoeff{m_n}{j}(1+nt)^{m_n-j}(2t)^j\mathbb{E}[|W_t|^j].
\end{align*}
Finally (arguing as for the proof of \eqref{eq:UB1MDP} above) we
take $\lambda=\frac{\tilde{x}}{2}t^{-\beta}$, and we have
\begin{align*}
\limsup_{t\to\infty}\frac{1}{t^{1-2\beta}}\log
&P\left(t^{\beta-1}(D_n(t)-(n-1)t)\geq
\tilde{x}\right)\\
\leq&\limsup_{t\to\infty}\frac{1}{t^{1-2\beta}}\left\{-\frac{\tilde{x}}{2}t^{-\beta}(\tilde{x}t^{1-\beta}+t(n-1))+
\kappa\left(\frac{\tilde{x}}{2}t^{-\beta}\right)t\right\}\\
=&\limsup_{t\to\infty}\left\{-\frac{\tilde{x}^2}{2}-\frac{\tilde{x}}{2}(n-1)t^\beta
+\frac{\tilde{x}}{2}t^\beta\left(\frac{\tilde{x}}{2}t^{-\beta}+n-1\right)\right\}=-\frac{\tilde{x}^2}{2}+\frac{\tilde{x}^2}{4}=-\frac{\tilde{x}^2}{4}.\
\Box
\end{align*}

\section{On the asymptotic behavior of some hitting probabilities}\label{sec:hitting-probabilities}
Let $\tau_{r_1}^{(n)}$ be the first hitting time of the
$n$-dimensional Euclidean Brownian motion on the sphere of radius
$r_1$ centered at the origin. Moreover, for $r>r_1$, let
$P_r(\tau_{r_1}^{(n)}<\infty)$ be the hitting probability when the
Euclidean Brownian motion starts at some point having Euclidean
distance $r$ from the origin. It is well known that
$$P_r(\tau_{r_1}^{(n)}<\infty)=\left\{\begin{array}{ll}
1&\ \mathrm{if}\ n=2\\
\frac{r^{2-n}}{r_1^{2-n}}&\ \mathrm{if}\ n\geq 3
\end{array}\right.\ (\mathrm{for}\ r>r_1).$$
Thus, as in immediate consequence, we have a \emph{polynomial
decay} as $r\to\infty$ in the fashion of large deviations, i.e.
$$\lim_{r\to\infty}\frac{1}{\log r}\log P_r(\tau_{r_1}^{(n)}<\infty)=-(n-2);$$
moreover the decay rate $w_e(n):=n-2$ is increasing with $n$ and
actually one expects that, the larger is the dimension of the
space, the faster is the decay of $P_r(\tau_{r_1}^{(n)}<\infty)$
as $r\to\infty$.

In this section we investigate the same kind of problem for the
hyperbolic Brownian motion. Let $T_{\eta_1}^{(n)}$ be the first
hitting time of the $n$-dimensional hyperbolic Brownian motion on
the hyperbolic sphere of radius $\eta_1$ centered at the origin
$O_n$. Moreover, for $\eta>\eta_1$, let
$P_\eta(T_{\eta_1}^{(n)}<\infty)$ be the hitting probability when
the hyperbolic Brownian motion starts at some point having
hyperbolic distance $\eta$ from $O_n$. It is known (see Corollary
3.1 and Corollary 3.2 in Cammarota and Orsingher
\cite{CammarotaOrsingher2012}) that, for $\eta>\eta_1$, we have
$$P_\eta(T_{\eta_1}^{(2)}<\infty)=\frac{\log\tanh\frac{\eta}{2}}{\log\tanh\frac{\eta_1}{2}},
\
P_\eta(T_{\eta_1}^{(3)}<\infty)=\frac{1-\coth\eta}{1-\coth\eta_1},$$
and, if we consider the values
$\{c(n,k):k\in\{0,\ldots,\frac{n-4}{2}\}\}$ defined by
$$c(n,0)=1\ \mbox{and}\ c(n,k)=\frac{(n-3)(n-5)\cdots (n-2k-1)}{(n-4)(n-6)\cdots (n-2k-2)},$$
we have
$$P_\eta(T_{\eta_1}^{(n)}<\infty)=\frac{\sum_{k=0}^{\frac{n-4}{2}}(-1)^kc(n,k)\frac{\cosh\eta}{\sinh^{n-2k-2}\eta}+
(-1)^{\frac{n}{2}}\frac{(n-3)!!}{(n-4)!!}\log\tanh\frac{\eta}{2}}{\sum_{k=0}^{\frac{n-4}{2}}(-1)^kc(n,k)\frac{\cosh\eta_1}{\sinh^{n-2k-2}\eta_1}+
(-1)^{\frac{n}{2}}\frac{(n-3)!!}{(n-4)!!}\log\tanh\frac{\eta_1}{2}}\
\mbox{for}\ n\in\{4,6,8,\ldots\}$$ and
$$P_\eta(T_{\eta_1}^{(n)}<\infty)=\frac{\sum_{k=0}^{\frac{n-5}{2}}(-1)^kc(n,k)\frac{\cosh\eta}{\sinh^{n-2k-2}\eta}+
(-1)^{\frac{n-5}{2}}\frac{(n-3)!!}{(n-4)!!}(1-\frac{\cosh\eta}{\sinh\eta})}{\sum_{k=0}^{\frac{n-5}{2}}(-1)^kc(n,k)\frac{\cosh\eta_1}{\sinh^{n-2k-2}\eta_1}+
(-1)^{\frac{n-5}{2}}\frac{(n-3)!!}{(n-4)!!}(1-\frac{\cosh\eta_1}{\sinh\eta_1})}\
\mbox{for}\ n\in\{5,7,9,\ldots\}.$$ An inspection for small values
of $n$ (some details will presented below, and one could have an
idea on how the methods can be adapted when $n$ is larger) lead us
to think that, for all $n\geq 2$, we have
\begin{equation}\label{eq:Lundberg-like-conjecture}
\lim_{\eta\to\infty}\frac{1}{\eta}\log
P_\eta(T_{\eta_1}^{(n)}<\infty)=-(n-1).
\end{equation}
Then we would have an \emph{exponential decay} as $\eta\to\infty$
in the fashion of large deviations and, as happens for
$w_e(n):=n-2$ above, the rate $w_h(n):=n-1$ is linearly increasing
with $n$. This kind of exponential decay has some analogy with
some asymptotic results in the literature for the logarithm of the
level crossing probabilities of real valued stochastic processes
with a \emph{stability condition}, which is satisfied when we have
a drift directed away from the level to cross (actually, as
pointed out in Remark \ref{rem:analogue-Euclidean-result}, the
hyperbolic Brownian motion has an implicit drift directed away
from the origin). Here we refer to the quite general result in
Duffy et al. \cite{DuffyLewisSullivan} (Theorem 2.2), but we have
in mind the simple case where $v$ and $a$ are the identity
function, $V=A=1$ and $J$ is convex (this situation comes up in
several cases; for the discrete time random walks with light tail
increments see Theorem 1 in Lehtonen and Nyrhinen
\cite{LehtonenNyrhinen}).\\

Now, we provide some details on the computations of the limit
\eqref{eq:Lundberg-like-conjecture} for $n\in\{2,\ldots,7\}$. We
start with the cases $n\in\{2,4,6\}$. When $n$ is even, we have to
handle the logarithmic term in the expression of
$P_\eta(T_{\eta_1}^{(n)}<\infty)$; in view of this, since
$\log\tanh\frac{\eta}{2}=-\log\frac{1+e^{-\eta}}{1-e^{-\eta}}$, by
eq. (1.513.1) in Gradshteyn and Ryzhik \cite{GradshteynRyzhik}
with $x=e^{-\eta}$, we get
$$\log\tanh\frac{\eta}{2}=-2\sum_{k=1}^\infty\frac{(e^{-\eta})^{2k-1}}{2k-1}.$$
For $n=2$ and $n=4$ we have to take into account
\begin{equation}\label{eq:bounds-for-series}
e^{-\eta}\leq\sum_{k=1}^\infty\frac{(e^{-\eta})^{2k-1}}{2k-1}\leq\frac{e^{-\eta}}{1-e^{-2\eta}}.
\end{equation}
Actually \eqref{eq:Lundberg-like-conjecture} holds with $n=2$ by
noting that
$$\liminf_{\eta\to\infty}\frac{1}{\eta}\log P_\eta(T_{\eta_1}^{(2)}<\infty)\geq\liminf_{\eta\to\infty}\frac{1}{\eta}\log(2e^{-\eta})=-1$$
and
$$\limsup_{\eta\to\infty}\frac{1}{\eta}\log P_\eta(T_{\eta_1}^{(2)}<\infty)\leq
\limsup_{\eta\to\infty}\frac{1}{\eta}\log\left(2\frac{e^{-\eta}}{1-e^{-2\eta}}\right)=-1;$$
moreover
$$\lim_{\eta\to\infty}\frac{1}{\eta}\log P_\eta(T_{\eta_1}^{(4)}<\infty)
=\lim_{\eta\to\infty}\frac{1}{\eta}\log\left(2\frac{e^\eta+e^{-\eta}}{(e^\eta-e^{-\eta})^2}-2\sum_{k=1}^\infty\frac{(e^{-\eta})^{2k-1}}{2k-1}\right)$$
and, again, we can prove \eqref{eq:Lundberg-like-conjecture} with
$n=4$ by deriving upper and lower bounds in terms of
\eqref{eq:bounds-for-series}. For $n=6$ we have
$$\lim_{\eta\to\infty}\frac{1}{\eta}\log P_\eta(T_{\eta_1}^{(6)}<\infty)
=\lim_{\eta\to\infty}\frac{1}{\eta}\log\left(8\frac{e^\eta+e^{-\eta}}{(e^\eta-e^{-\eta})^4}-3\frac{e^\eta+e^{-\eta}}{(e^\eta-e^{-\eta})^2}
+3\sum_{k=1}^\infty\frac{(e^{-\eta})^{2k-1}}{2k-1}\right)$$ and we
are able to prove \eqref{eq:Lundberg-like-conjecture} by deriving
upper and lower bounds as in the previous cases, by means of the
following relationships
$$e^{-\eta}+\frac{e^{-3\eta}}{3}\leq\sum_{k=1}^\infty\frac{(e^{-\eta})^{2k-1}}{2k-1}\leq e^{-\eta}+\frac{e^{-3\eta}}{3}+\frac{e^{-5\eta}}{1-e^{-2\eta}}.$$
instead of \eqref{eq:bounds-for-series}. Finally the cases
$n\in\{3,5,7\}$:
\begin{align*}
\lim_{\eta\to\infty}\frac{1}{\eta}\log
P_\eta(T_{\eta_1}^{(3)}<\infty)=&\lim_{\eta\to\infty}\frac{1}{\eta}\log(\coth\eta-1)\\
=&\lim_{\eta\to\infty}\frac{1}{\eta}\log\left(\frac{e^\eta+e^{-\eta}}{e^\eta-e^{-\eta}}-1\right)
=\lim_{\eta\to\infty}\frac{1}{\eta}\log\left(\frac{2}{e^{2\eta}-1}\right)=-2;
\end{align*}
\begin{align*}
\lim_{\eta\to\infty}\frac{1}{\eta}\log
P_\eta(T_{\eta_1}^{(5)}<\infty)=&\lim_{\eta\to\infty}\frac{1}{\eta}\log\left(\frac{\cosh\eta}{\sinh^3\eta}+2\left(1-\frac{\cosh\eta}{\sinh\eta}\right)\right)\\
=&\lim_{\eta\to\infty}\frac{1}{\eta}\log\left(4\frac{e^\eta+e^{-\eta}}{(e^\eta-e^{-\eta})^3}-4\frac{e^{-\eta}}{e^\eta-e^{-\eta}}\right)\\
=&\lim_{\eta\to\infty}\frac{1}{\eta}\log\left(\frac{e^{-\eta}+o(e^{-\eta})}{(e^\eta-e^{-\eta})^3}\right)=-4;
\end{align*}
\begin{align*}
\lim_{\eta\to\infty}\frac{1}{\eta}\log
P_\eta(T_{\eta_1}^{(7)}<\infty)=&\lim_{\eta\to\infty}\frac{1}{\eta}\log\left(\frac{\cosh\eta}{\sinh^5\eta}-
\frac{4}{3}\frac{\cosh\eta}{\sinh^3\eta}-\frac{8}{3}\left(1-\frac{\cosh\eta}{\sinh\eta}\right)\right)\\
=&\lim_{\eta\to\infty}\frac{1}{\eta}\log\left(16\frac{e^\eta+e^{-\eta}}{(e^\eta-e^{-\eta})^5}
-\frac{16}{3}\frac{e^\eta+e^{-\eta}}{(e^\eta-e^{-\eta})^3}+\frac{16}{3}\frac{e^{-\eta}}{e^\eta-e^{-\eta}}\right)\\
=&\lim_{\eta\to\infty}\frac{1}{\eta}\log\left(\frac{e^{-\eta}+o(e^{-\eta})}{(e^\eta-e^{-\eta})^5}\right)=-6.
\end{align*}

\section{On a recent LDP in the literature}\label{sec:Hirao}
Theorem 1.1 in Hirao \cite{Hirao} provides the LDP for
$\left\{\frac{d_{\mathbb{F}}\left(Z_z^{\mathbb{F}}(t),z\right)}{t}:t>0\right\}$,
where $\left\{Z_z^{\mathbb{F}}(t):t\geq 0\right\}$ is a Brownian
motion on a hyperbolic space $H^n(\mathbb{F})$ with distance
function $d_{\mathbb{F}}$ (several choices of $\mathbb{F}$ are
allowed), and $Z_z^{\mathbb{F}}(0)=z$. Here we want to illustrate
the relationship between Theorem 1.1 in Hirao \cite{Hirao} when
$\mathbb{F}$ concerns the parameters $(k,m)=(0,n-1)$ in eq. (1) in
Hirao \cite{Hirao}, and Proposition \ref{prop:LDP} in this paper.
Actually, if we denote the rate function for the first LDP by
$\Lambda^*$, we have $\Lambda^*(x)=2I_1(x)$ for all $x\geq 0$.

Firstly we should have $\Lambda^*(x)=I_1(x)$ for all $x\geq 0$ if
we consider $\frac{t}{2}$ in place of $t$ because the exponential
part of function $h_{\mathbb{F}}^{(k,m)}$ in eq. (4) in Hirao
\cite{Hirao} with $(k,m)=(0,n-1)$ would coincide with the one of
$h_n(\eta,\frac{t}{2})$ in this paper. Furthermore the proof of
Theorem 1.1 in Hirao \cite{Hirao} shows the existence of the limit
\begin{equation}\label{eq:Hirao-limit}
\lim_{t\to\infty}\frac{1}{t}\log\mathbb{E}\left[e^{\lambda
d_{\mathbb{F}}\left(Z_z^{\mathbb{F}}(t),z\right)}\right]=\frac{1}{2}\lambda(\lambda+2k+m)=:\Lambda(\lambda)
\end{equation}
for all $\lambda\in\mathbb{R}$ and, by an application of the
G\"{a}rtner Ellis Theorem, the LDP holds with the good rate
function $\Lambda^*$ defined by
\begin{equation}\label{eq:def-Lambda-conjugate}
\Lambda^*(x):=\sup_{\lambda\in\mathbb{R}}\{\lambda
x-\Lambda(\lambda)\}.
\end{equation}
Thus, by \eqref{eq:Hirao-limit} and
\eqref{eq:def-Lambda-conjugate}, the rate function in Theorem 1.1
in Hirao \cite{Hirao} is
$$\Lambda^*(x)=\frac{1}{2}\left(x-\frac{2k+m}{2}\right)^2$$
for all $x\in\mathbb{R}$.

We remark that, since the random variables
$\left\{d_{\mathbb{F}}\left(Z_z^{\mathbb{F}}(t),z\right):t>0\right\}$
should be nonnegative, there is a slight inexactness in this
proof. Actually, if a family of nonnegative random variables
$\{Z(t):t>0\}$ satisfies the LDP with a rate function $I$, the
lower bound for the open set $G=(-\infty,0)$ in the definition of
LDP would yield $I(x)=\infty$ for $x\in(-\infty,0)$; moreover, for
the same reason, the limit \eqref{eq:Hirao-limit} fails because
the function $\Lambda$ should be nondecreasing. We think that the
proof could be corrected showing that we have
\begin{equation}\label{eq:Hirao-limit-correction}
\lim_{t\to\infty}\frac{1}{t}\log\mathbb{E}\left[e^{\lambda
d_{\mathbb{F}}\left(Z_z^{\mathbb{F}}(t),z\right)}\right]=
\left\{\begin{array}{ll}
\frac{1}{2}\lambda(\lambda+2k+m)&\ \mbox{if}\ \lambda\geq -\frac{2k+m}{2}\\
-\frac{1}{2}\left(\frac{2k+m}{2}\right)^2&\ \mbox{if}\
\lambda<-\frac{2k+m}{2}
\end{array}\right.=:\Lambda(\lambda)
\end{equation}
in place of \eqref{eq:Hirao-limit}; actually, in such a case, the
hypotheses of G\"{a}rtner Ellis Theorem would be satisfied and the
LDP would hold with rate function $\Lambda^*$ defined by
\eqref{eq:def-Lambda-conjugate}. Moreover the function $\Lambda$
in \eqref{eq:Hirao-limit-correction} is nondecreasing, and
\eqref{eq:def-Lambda-conjugate} and
\eqref{eq:Hirao-limit-correction} yield
$$\Lambda^*(x)=\left\{\begin{array}{ll}
\frac{1}{2}\left(x-\frac{2k+m}{2}\right)^2&\ \mathrm{if}\ x\geq 0\\
\infty&\ \mathrm{if}\ x<0
\end{array}\right.$$
(thus $\Lambda^*(x)=\infty$ for $x\in(-\infty,0)$).

\end{document}